 \documentclass[12pt]{article}
 \usepackage{amsmath,amsthm,amssymb,amscd}
\title{Mukai implies McKay:\\
 the McKay correspondence as an\\ equivalence of derived categories}
\author{Tom Bridgeland \and Alastair King \and Miles Reid}
\date{To Andrei Tyurin on his 60th birthday}
\newcommand{\Step}[1]{\paragraph{\sc{Step #1}}}

\newtheorem{thm}{Theorem}[section]
\newtheorem{cor}[thm]{Corollary}
\newtheorem{lemma}[thm]{Lemma}
\newenvironment{pf}{\paragraph{Proof}}{\qed\par\medskip}

 \newcommand{\1}{^{-1}}
 \newcommand{\iso}{\cong}
 \newcommand{\comp}{\circ}
 \newcommand{\blank}{\kern.0mm\smash{-}\kern.0mm}
 \newcommand{\blob}{{\scriptscriptstyle\bullet}}
 \newcommand{\dual}{\vee}
 \newcommand{\tensor}{\otimes}

 \newcommand{\FF}{\mathbb F}
 \newcommand{\QQ}{\mathbb Q}
 \newcommand{\ZZ}{\mathbb Z}
 \newcommand{\C}{\mathbb C}
 \newcommand{\Z}{\mathbb Z}

 \newcommand{\K}{\mathcal K}
 \newcommand{\A}{\mathcal A}
 \newcommand{\B}{\mathcal B}
 \newcommand{\sP}{\mathcal P}
 \newcommand{\sR}{\mathcal R}
 \newcommand{\sS}{\mathcal S}
 \newcommand{\sZ}{\mathcal Z}
 \newcommand{\Oh}{\mathcal O}
 \newcommand{\sQ}{\operatorname{\mathcal Q}}

 \newcommand{\ep}{\varepsilon}
 \newcommand{\fie}{\varphi}
 \newcommand{\om}{\omega}
 \newcommand{\Om}{\Omega}
 \newcommand{\Ups}{\Upsilon}

 \newcommand{\Ybar}{\overline{Y}}
 \newcommand{\m}{m}
 \newcommand{\Se}{S}
 \newcommand{\lift}[2]{\lambda_{#1}^{#2}}
 \newcommand{\Phic}{\Phi_{\mathrm c}}

 \newcommand{\qc}{\mathrm{qc}}
 
 \newcommand{\id}{\operatorname{id}}
 \newcommand{\Qco}{\operatorname{Qcoh}}
 \newcommand{\Aut}{\operatorname{Aut}} 
 \newcommand{\Coh}{\operatorname{Coh}}
 \newcommand{\D}{\operatorname{D}}
 \newcommand{\Dc}{\operatorname{D_{\mathrm c}}}
 \newcommand{\DcG}{\operatorname{D_{\mathrm c}^G}}
 \newcommand{\GExt}{\operatorname{\hbox{$G$}-Ext}}
 \newcommand{\GHom}{\operatorname{\hbox{$G$}-Hom}}
 \newcommand{\Ext}{\operatorname{Ext}}
 \newcommand{\Hom}{\operatorname{Hom}}
 \newcommand{\sHom}{\operatorname{\mathcal H\mathit{om}}}
 \newcommand{\AHilb}{\operatorname{\hbox{$A$}-Hilb}}
 \newcommand{\GHilb}{\operatorname{\hbox{$G$}-Hilb}}
 \newcommand{\Ob}{}
 \newcommand{\Supp}{\operatorname{Supp}}
 \newcommand{\Spec}{\operatorname{Spec}}
 \newcommand{\SL}{\operatorname{SL}}
 
 \newcommand{\eu}{\operatorname{\chi}}
 \newcommand{\hd}{\operatorname{hom\, dim}}
 \newcommand{\codim}{\operatorname{codim}}
 \renewcommand{\L}{\mathbf L}
 \newcommand{\R}{\mathbf R}
 \newcommand{\Ltensor}{\stackrel{\mathbf L}{\tensor}}
 \newcommand{\coker}{\operatorname{coker}}
 \newcommand{\Ker}{\operatorname{Ker}}
 \renewcommand{\Im}{\operatorname{Im}}

 \newcommand{\into}{\hookrightarrow}
 \newcommand{\lRa}[1]{\xrightarrow{\ #1\ }}
 \newcommand{\lra}{\longrightarrow}
 
 \newcommand{\bua}{\bigg{\uparrow}}

\makeatletter
 \newlength{\typesize}
\makeatother

\newlength{\vvoff}
\newlength{\hhoff}

\newcommand{\locateoffcenter}[1]{%
\addtolength{\vvoff}{-0.25\typesize}%
\raisebox{\vvoff}{\hspace{\hhoff}\makebox(0,0){\smash{#1}}}
}
\newcommand{\object}[1]{%
\setlength{\vvoff}{0pt}%
\setlength{\hhoff}{0pt}%
\locateoffcenter{#1}
}

\newcommand{\swlabel}[1]{%
\setlength{\vvoff}{-0.5\typesize}%
\setlength{\hhoff}{0.75\typesize}%
\locateoffcenter{#1}
}
\newcommand{\nwlabel}[1]{%
\setlength{\vvoff}{-0.5\typesize}%
\setlength{\hhoff}{-0.75\typesize}%
\locateoffcenter{#1}
}
\newcommand{\selabel}[1]{%
\setlength{\vvoff}{0.5\typesize}%
\setlength{\hhoff}{0.75\typesize}%
\locateoffcenter{#1}
}
\newcommand{\nelabel}[1]{%
\setlength{\vvoff}{0.5\typesize}%
\setlength{\hhoff}{-0.75\typesize}%
\locateoffcenter{#1}
}
\begin{document}
\maketitle

\begin{abstract}
Let $G$ be a finite group of automorphisms of a non\-singular complex
threefold $M$ such that the canonical bundle $\om_M$ is locally
trivial as a $G$-sheaf. We prove that the Hilbert scheme
$Y=\GHilb{M}$ parametrising
$G$-clusters in $M$ is a crepant resolution of $X=M/G$ and that there
is a derived equi\-valence (Fourier--Mukai transform) between coherent
sheaves on $Y$ and coherent \hbox{$G$-sheaves} on $M$. This identifies
the K~theory of $Y$ with the equi\-variant K~theory of $M$, and thus
generalises the classical McKay correspondence. Some higher dimensional
extensions are possible.
\end{abstract}

%=======================================================================
\section{Introduction}
%=======================================================================

The classical McKay correspondence relates representations of a finite
subgroup $G\subset\SL(2,\C)$ to the cohomology of the well-known minimal
resolution of the Kleinian singularity $\C^2/G$. Gonzalez-Sprinberg and
Verdier \cite{GV} interpreted the McKay correspondence as an isomorphism
on K~theory, observing that the representation ring of $G$ is equal to
the $G$-equi\-variant K~theory of $\C^2$.

A natural generalisation is to replace $\C^2$ by a non\-singular
quasi\-projective complex variety $M$ of dimension $n$ and $G$ by
a finite group
of automorphisms of $M$, with the property that the stabiliser subgroup of
any point $x\in M$ acts on the tangent space $T_xM$ as a subgroup of
$\SL(T_xM)$. Thus the canonical bundle $\om_M$ is locally trivial as a
$G$-sheaf, in the sense that every point of $M$ has a $G$-invariant
open neighbourhood on which there is a nonvanishing $G$-invariant
$n$-form. This implies that the quotient variety $X=M/G$ has only
Gorenstein singularities.

The natural generalisation of the McKay correspondence should then be an
isomorphism between the $G$-equi\-variant K~theory of $M$ and the ordinary 
K~theory of a crepant resolution $Y$ of $X$, that is, a resolution of
singularities $\tau\colon Y\to X$ such that $\tau^*(\om_X)=\om_Y$. Crepant
resolutions of Gorenstein quotient singularities are known to exist in
dimension $n=3$, but only through a case by case analysis of the local
linear actions by Ito, Markushevich and Roan (see Roan \cite{Roa} and
references given there). In dimension $\ge4$, crepant resolutions exist
only in rather special cases.

The point of view of this paper is that the derived category is the
natural context for this formulation of the correspondence, and,
more importantly, provides key tools for an appropriately general proof.
Indeed, this point of view is not so revolutionary. 
Gonzalez-Sprinberg and Verdier were aware that their isomorphism
on K~theory would lift to a derived equivalence and an explicit proof
of this was given by Kapranov and Vasserot \cite{KV}.
Moreover, the statement of the McKay correspondence in 3 dimensions in terms 
of K~theory and derived categories is contained in 
Reid \cite[Conjecture~4.1]{R}.
One surprise, however, is that the
methods of the derived category are powerful enough to prove the
existence of a crepant resolution in 3 dimensions, without any case
by case analysis.

\smallskip

A good candidate for a crepant resolution of $X$ is Nakamura's $G$-Hilbert
scheme $\GHilb{M}$ parametrising $G$-clusters or `scheme theoretic
$G$-orbits' on $M$: recall that a {\em cluster} $Z\subset M$ is a zero
dimensional subscheme, and a {\em $G$-cluster} is a $G$-invariant cluster
whose global sections $H^0(\Oh_Z)$ are isomorphic to the regular 
representation $\C[G]$ of $G$. Clearly, a $G$-cluster has length $|G|$ and
a free $G$-orbit is a $G$-cluster. There is a Hilbert--Chow morphism
 \[
 \tau\colon\GHilb{M}\lra X,
 \]
which, on closed points, sends a $G$-cluster to the orbit supporting it.
Note that $\tau$~is a projective morphism, is onto and is birational on one
component.

When $M=\C^3$ and $G\subset\SL(3,\C)$ is Abelian, Nakamura \cite{N} proved
that $\GHilb{M}$ is irreducible and is a crepant resolution of $X$ (compare
also Reid \cite{R} and Craw and Reid \cite{CR}). He conjectured that the same
result holds for an arbitrary finite subgroup $G\subset\SL(3,\C)$. 
Ito and Nakajima \cite{IN} observed that the 
construction of Gonzalez-Sprinberg and Verdier \cite{GV} 
is the $M=\C^2$ case of a natural correspondence between the
equi\-variant K~theory of $M$ and the ordinary K~theory of $\GHilb{M}$. 
They proved that this correspondence is an isomorphism when $M=\C^3$ and
$G\subset\SL(3,\C)$ is Abelian by constructing an explicit resolution
of the diagonal in Beilinson style. Our approach via Fourier--Mukai
transforms leaves this resolution of the diagonal implicit (it appears
as the object $\sQ$ of $\D(Y\times Y)$ in
Section~\ref{Sec6!Proj_case}), and seems to give a more direct argument.
Two of the main consequences of the results of this
paper are that Nakamura's conjecture is true and that the natural
correspondence on K~theory is an isomorphism for all finite subgroups of
$\SL(3,\C)$.

As already indicated, the basic approach of the paper is to lift the McKay
correspondence to the appropriate derived categories. We may then apply the
techniques of Fourier--Mukai transforms, in particular the ideas of
Bridgeland \cite{Br1} and \cite{Br2}, to show that it is an equivalence at
this level. The more formal nature of the arguments means that they work
equally well for arbitrary quasi\-projective varieties. In fact, they are
somewhat simpler for projective varieties, and we therefore deal with this
case first.

Since it is not known whether $\GHilb{M}$ is irreducible or even connected
in general, we take as our initial candidate for a resolution $Y$ the
{\em irreducible component} of $\GHilb{M}$ containing the free $G$-orbits,
that is, the component mapping birationally to $X$. The aim is to show that
$Y$ is a crepant resolution, and to construct an equivalence between the
derived categories $\D(Y)$ of coherent sheaves on $Y$ and $\D^G(M)$ of
coherent \hbox{$G$-sheaves} on $M$. A particular consequence of
this equivalence is that $Y=\GHilb{M}$ when $M$ has dimension 3.

\medskip

We now describe the correspondence and our results in more detail.
Let $M$ be a non\-singular quasiprojective complex variety of dimension $n$ and
let $G\subset\Aut(M)$ be a finite group of automorphisms of $M$
such that $\omega_M$
is locally trivial as a $G$-sheaf. Put $X=M/G$ and let
$Y\subset\GHilb{M}$ be the
irreducible component containing the free orbits, as described above.
Write $\sZ$ for the universal closed subscheme $\sZ\subset Y\times M$
and $p$ and $q$ for its projections to $Y$ and $M$. There is a commutative
diagram of schemes
\begin{equation*}\label{maindiag}
 \setlength{\unitlength}{36pt}
 \begin{picture}(2.2,2.2)(0,0)
 \put(1,0){\object{$X$}}
 \put(0,1){\object{$Y$}}
 \put(2,1){\object{$M$}}
 \put(1,2){\object{$\sZ$}}
 \put(0.25,0.75){\vector(1,-1){0.5}}
 \put(0.5,0.5){\nwlabel{$\tau$}}
 \put(1.75,0.75){\vector(-1,-1){0.5}}
 \put(1.5,0.5){\swlabel{$\pi$}}
 \put(0.75,1.75){\vector(-1,-1){0.5}}
 \put(0.5,1.5){\nelabel{$p$}}
 \put(1.25,1.75){\vector(1,-1){0.5}}
 \put(1.5,1.5){\selabel{$q$}}
 \end{picture}
\end{equation*}
in which $q$ and $\tau$ are birational, $p$ and $\pi$ are finite, and $p$ is
flat. Let $G$ act trivially on $Y$ and $X$, so that all morphisms
in the diagram are equi\-variant.

Define the functor
 \[
 \Phi=\R q_*\comp p^*\colon\D(Y)\lra \D^G(M),
 \]
where a sheaf $E$ on $Y$ is viewed as a $G$-sheaf by giving it the trivial
action. Note that $p^*$ is already exact, so we do not need to write
$\L p^*$. Our main result is the following.

\medskip\goodbreak

\begin{thm}
\label{second}
Suppose that the fibre product
 \[
 Y\times_X Y= \Bigl\{(y_1, y_2)\in Y\times Y \Bigm|
 \tau(y_1)=\tau(y_2)\Bigr\} \subset Y\times Y
\]
has dimension\/ $\le n+1$. 
Then\/ $Y$ is a crepant resolution of\/ $X$ and\/ $\Phi$ is an
equivalence of categories.
\end{thm}

When $n\le 3$ the condition of the theorem always holds because the
exceptional locus of $Y\to X$ has dimension $\le2$. In this case we
can also show that $\GHilb{M}$ is irreducible, so we obtain

\begin{thm}
\label{IN}
Suppose $n\le3$. Then $\GHilb{M}$ is irreducible and is a crepant
resolution of $X$, and\/ $\Phi$ is an equivalence of categories.
\end{thm}

The condition of Theorem~\ref{second} also holds whenever $G$
preserves a complex symplectic form on $M$ and $Y$ is a crepant
resolution of $X$, because such a resolution is symplectic and hence
semismall (see
Verbitsky \cite{Vb}, Theorem~2.8 and compare Kaledin \cite{Ka}).

\begin{cor}
\label{Kal}
Suppose $M$ is a complex symplectic variety
and $G$ acts by symplectic automorphisms.
Assume that $Y$ is a crepant resolution of $X$.
Then $\Phi$ is an equivalence of categories.
\end{cor}

Note that the condition of Theorem~\ref{second} certainly fails in
dimension $\ge4$ whenever $Y\to X$ has an exceptional divisor over a
point. This is to be expected since there are many examples of finite
subgroups $G\subset\SL(4,\C)$ for which the corresponding quotient
singularity $\C^4/G$ has no crepant resolution.

\subsection*{Acknowledgements}
The first author would like to thank the ICTP, Trieste and EPSRC
for financial support whilst this paper was written.

%=======================================================================
\section{Category theory}
%=======================================================================

This section contains some basic category theory, most of which is well
known. The only non\-trivial part is Section~\ref{trings} where we state a
condition for an exact functor between triangulated categories to be an
equivalence.
 
\subsection{Triangulated categories}

A triangulated category is an additive category $\A$ equipped with a
{\em shift auto\-morphism} $T_{\A}\colon\A\to\A\colon a\mapsto a[1]$
and a collection of {\em distinguished triangles}
 \[
 a_1\lRa{f_1}a_2\lRa{f_2}a_3\lRa{f_3}a_1[1]
 \]
of morphisms of $\A$ satisfying certain axioms (see Verdier \cite{V}).
We write $a[i]$ for $T_{\A}^i(a)$ and
 \[
 \Hom^i_{\A}(a_1,a_2)=\Hom_{\A}(a_1,a_2[i]).
 \]
A triangulated category $\A$ is {\em trivial}\/ if every object
is a zero object.

The principal example of a triangulated category is the derived category
$\D(A)$ of an Abelian category $A$. An object of $\D(A)$ is a bounded
complex of objects of $A$ up to quasi-isomorphism, the shift functor
moves a complex to the left by one place and a distinguished triangle is
the mapping cone of a morphism of complexes. 
In this case, for objects $a_1,a_2\in A$, %%%
one has $\Hom^i_{\D(A)}(a_1,a_2)=\Ext^i_A(a_1,a_2)$.

A functor $F\colon\A\to\B$ between triangulated categories is {\em exact}\/
if it commutes with the shift automorphisms and takes distinguished
triangles of $\A$ to distinguished triangles of $\B$. For example, derived
functors between derived categories are exact.

\subsection{Adjoint functors}

Let $F\colon\A\to\B$ and $G\colon\B\to\A$ be functors.
An adjunction for $(G,F)$ is a bifunctorial isomorphism
 \[
 \Hom_{\A}(G\blank,\blank)\,\iso\,\Hom_{\B}(\blank,F\blank).
 \]
In this case, we say that $G$ is left adjoint to $F$ or that $F$ is right
adjoint to $G$. When it exists, a left or right adjoint to a given functor
is unique up to isomorphism of functors. The adjoint of a composite functor
is the composite of the adjoints. An adjunction determines and is determined
by two natural transformations $\ep\colon G\comp F\to\id_\A$ and
$\eta\colon\id_\B\to F\comp G$ that come from applying the adjunction to
$1_{Fa}$ and $1_{Gb}$ respectively (see Mac~Lane \cite[IV.1]{Mac} for more
details).

The basic adjunctions we use in this paper are described in
Section~\ref{adjunctions} below.

\subsection{Fully faithful functors and equivalences}
\label{fff}

A functor $F\colon\A\to \B$ is {\em fully faithful}\/ if for any pair of
objects $a_1$, $a_2$ of $\A$, the map
 \[
 F\colon\Hom_{\A}(a_1,a_2)\to\Hom_{\B}(Fa_1,Fa_2)
 \]
 is an isomorphism. One should think of $F$ as an `injective' functor.
This is more clear when $F$ has a left adjoint $G\colon\B\to\A$
(or a right adjoint $H\colon\B\to\A$), 
in which case $F$ is fully faithful if and only if the
natural transformation $G\comp F\to \id_{\A}$ (or $\id_{\A}\to H\comp F$)
is an isomorphism.

A functor $F$ is an {\em equivalence} if there is an `inverse' functor
$G\colon\B\to\A$ such that $G\comp F\iso \id_{\A}$ and $F\comp G\iso
\id_{\B}$. In this case $G$ is both a left and right adjoint to $F$ (see
Mac~Lane \cite[IV.4]{Mac}). In practice, we show that $F$ is an equivalence
by writing down an adjoint (a priori, one-sided) and proving that it is an
inverse.
One simple example of this is the following.

\begin{lemma}
\label{extra}
Let $\A$ and $\B$ be triangulated categories and $F\colon\A\to \B$ a
fully faithful exact functor with a right adjoint $H\colon\B\to\A$.
Then $F$ is an equivalence 
if and only if $Hc\iso 0$ implies $c\iso 0$ for any object $c\in\Ob\B$.
\end{lemma}

\begin{pf}
By assumption $\eta\colon\id_\A\to H\comp F$ is an
isomorphism, so $F$ is an equivalence if and only if
$\ep\colon F\comp H\to\id_\B$ is an isomorphism.
Thus the `only if' part of the lemma is immediate, since $c\iso FHc$.

For the `if' part, take any object $b\in\Ob\B$ and
embed the natural adjunction map $\ep_b$ in a triangle
\begin{equation}
\label{semiortho}
 c\to FHb\lRa{\ep_b} b\to c[1].
\end{equation}
If we apply $H$ to this triangle, then $H(\ep_b)$ is an isomorphism,
because $\eta_{Hb}$ is an isomorphism and
$H(\ep_b)\comp\eta_{Hb}=1_{Hb}$ (\cite[IV.1, Theorem 1]{Mac}).
Hence $Hc\iso 0$ and so $c\iso 0$ by hypothesis. 
Thus $\ep_b$ is an isomorphism, as required.
\end{pf}

One may understand this lemma in a broader context
as follows.
The triangle~(\ref{semiortho}) shows that,
when $F$ is fully faithful with right adjoint $H$,
there is a `semi-orthogonal' decomposition $\B=(\Im F,\Ker H)$,
where 
\begin{align*}
  \Im F &= \{ b\in\Ob\B : \text{$b\iso Fa$ for some $a\in\Ob\A$} \},\\  
  \Ker H &= \{ c\in\Ob\B : Hc\iso 0 \}.
\end{align*}
Since $F$ is fully faithful, 
the fact that $b\iso Fa$ for some object $a\in\Ob\A$
necessarily means that $b\iso FHb$,
so only zero objects are in both subcategories.
The semi-orthogonality condition also requires that 
$\Hom_{\B}(b,c)=0$ for all $b\in\Im F$ and $c\in\Ker H$, %%%
which is immediate from the adjunction.
The lemma then has the very reasonable interpretation that
if $\Ker H$ is trivial,
then $\Im F=\B$ and $F$ is an equivalence.
Note that if $G$ is a left adjoint for $F$, then there is a similar
semi-orthogonal decomposition on the other side 
$\B=(\Ker G,\Im F)$ and a corresponding version of the lemma.
For more details on semi-orthogonal decompositions see Bondal \cite{Bo}.

\subsection{Spanning classes and orthogonal decomposition}

A {\em spanning class} for a triangulated category $\A$ is a subclass
$\Om$ of the objects of $\A$ such that for any object $a\in\Ob\A$
 \[
 \Hom^i_{\A}(a,\om)=0 \quad\mbox{for all }
 \om\in\Om, i\in\Z\quad\mbox{implies }a\iso 0
 \]
and
 \[
 \Hom^i_{\A}(\om,a)=0\quad\mbox{for all }
 \om\in\Om, i\in\Z\quad\mbox{implies }a\iso 0.
 \]
For example, the set of skyscraper sheaves $\{\Oh_x\colon x\in X\}$
on a non\-singular variety $X$ is a spanning class for $\D(X)$.

A triangulated category $\A$ is {\em decomposable} as an orthogonal direct
sum of two full subcategories $\A_1$ and $\A_2$ if every object of $\A$ is
isomorphic to a direct sum $a_1\oplus a_2$ with $a_j\in\Ob\A_j$, and if
 \[
 \Hom_{\A}^i(a_1,a_2)=\Hom_{\A}^i(a_2,a_1)=0
 \]
for any pair of objects $a_j\in\Ob\A_j$ and all integers $i$.
The category $\A$ is indecomposable if for any such decomposition one
of the two subcategories $\A_i$ is trivial.
For example, if $X$ is a scheme, $\D(X)$ is indecomposable
precisely when $X$ is connected. For more details see Bridgeland \cite{Br1}.

\subsection{Serre functors}

The properties of Serre duality on a non\-singular projective variety
were abstracted by Bondal and Kapranov \cite{BK} into the 
notion of a Serre functor on a triangulated category.
Let $\A$ be a triangulated category in which all the $\Hom$ sets are
finite dimensional vector spaces. A {\em Serre functor} for $\A$ is an
exact equivalence $\Se\colon\A\to\A$ inducing bifunctorial
iso\-morphisms
 \[
 \Hom_{\A}(a,b)\to\Hom_{\A}(b,\Se (a))^{\dual}
 \quad \text{for all $a,b\in\Ob\A$}
 \]
that satisfy a simple compatibility condition (see \cite{BK}). When a
Serre functor exists, it is unique up to isomorphism of functors. We say
that $\A$ has {\em trivial}\/ Serre functor if for some integer $i$ the
shift functor $[i]$ is a Serre functor for $\A$.

The main example is the bounded derived category of coherent sheaves
$\D(X)$ on a non\-singular projective $n$-fold $X$, having the Serre
functor
 \[
 S_X(\blank)=(\blank\tensor\om_X)[n].
 \]
Thus $\D(X)$ has trivial Serre functor if and only if the canonical bundle
of $X$ is trivial.

\subsection{A criterion for equivalence}
\label{trings}

Let $F\colon\A\to\B$ be an exact functor between triangulated categories
with Serre functors $\Se_{\A}$ and $\Se_{\B}$. Assume that $F$ has a left
adjoint $G\colon\B\to\A$. Then $F$ also has a right adjoint 
$H=\Se_{\A}\comp G\comp \Se_{\B}^{-1}$.

\begin{thm}
\label{tring}
Suppose there is a spanning class $\Om$ for $\A$ such that
\begin{equation*}
F\colon\Hom^i_{\A}(\om_1,\om_2)\to \Hom^i_{\B}(F\om_1,F\om_2)
\end{equation*}
is an isomorphism for all $i\in\Z$ and all
$\om_1,\om_2\in\Om$. 
Then $F$ is fully faithful.
\end{thm}

\begin{pf}
See \cite[Theorem 2.3]{Br1}.
\end{pf}

\begin{thm}
\label{tring2}
Suppose further that $\A$ is non\-trivial, 
that $\B$ is indecomposable and
that $F\Se_{\A}(\om)\iso\Se_{\B}F(\om)$ for all $\om\in\Om$.
Then $F$ is an equivalence of categories.
\end{thm}

\begin{pf}
Consider an object $b\in\Ob\B$. 
For any $\om\in\Om$ and $i\in\Z$ we have isomorphisms
\begin{eqnarray*}
 &\Hom_{\A}^i(\om,Gb)=\Hom_{\A}^i(Gb,\Se_{\A} \om)^{\dual}=
\Hom_{\B}^i(b,F\Se_{\A} \om)^{\dual} \\
&=\Hom_{\B}^i(b,\Se_{\B} F \om)^{\dual}=
\Hom_{\B}^i(F\om,b)=\Hom_{\A}^i(\om,Hb),
\end{eqnarray*}
using Serre duality and the adjunctions for $(G,F)$ and $(F,H)$. Since
$\Om$ is a spanning class we can conclude that $Gb\iso 0$ precisely when
$Hb\iso 0$. Then the result follows from \cite[Theorem 3.3]{Br1}.
\end{pf}

The proof of Theorem 3.3 in \cite{Br1} may
be understood as follows.
If $\Ker H\subset\Ker G$, then 
the semiorthogonal decomposition described 
at the end of Section~\ref{fff} becomes an orthogonal decomposition.
Hence $\Ker H$ must be trivial, because $\B$ is indecomposable and
$\A$, and hence $\Im F$, is nontrivial.
Thus $\Im F=\B$ and $F$ is an equivalence.

%=======================================================================
\section{Derived categories of sheaves}
%=======================================================================

This section is concerned with various general properties of complexes
of $\Oh_X$-modules on a scheme $X$. Note that all our schemes are of
finite type over $\C$. Given a scheme $X$, define $\D^{\qc}(X)$ to be
the (unbounded) derived category of the Abelian category 
$\Qco(X)$ of quasi\-coherent sheaves on $X$. 
Also define $\D(X)$ to be the full subcategory of
$\D^{\qc}(X)$ consisting of complexes with bounded and coherent cohomology.

\subsection{Geometric adjunctions}
\label{adjunctions}
Here we describe three standard adjunctions that arise in algebraic
geometry and are used frequently in what follows. For the first example,
let $X$ be a scheme and $E\in\D(X)$ an object of finite homological
dimension. Then the derived dual
 \[
 E^{\dual}=\R\sHom_{\Oh_X}(E,\Oh_X)
 \]
also has finite homological dimension, and the functor $\blank\Ltensor E$ is
both left and right adjoint to the functor $\blank\Ltensor E^{\dual}$.

For the second example take a morphism of schemes $f\colon X\to Y$. The
functor
 \[
 \R f_*\colon \D^{\qc}(X)\lra \D^{\qc}(Y)
 \]
has the left adjoint
 \[
 \L f^*\colon \D^{\qc}(Y)\lra \D^{\qc}(X).
 \]
If $f$ is proper then $\R f_*$ takes $\D(X)$ into $\D(Y)$. If $f$
has finite Tor dimension (for example if $f$ is flat, or $Y$ is
non\-singular) then $\L f^*$ takes $\D(Y)$ into $\D(X)$.

The third example is Grothendieck duality. Again take a morphism of
schemes $f\colon X\to Y$. The functor $\R f_*$ has a right adjoint
 \[
 f^!\colon \D^{\qc}(Y)\lra \D^{\qc}(X)
 \]
and moreover, if $f$ is proper and of finite Tor dimension, there is an
isomorphism of functors
\begin{equation}
\label{amos}
f^!(\blank)\,\iso\, \L f^*(\blank)\Ltensor f^!(\Oh_Y).
\end{equation}
Neeman \cite{Ne} has recently given a completely formal proof of these
statements in terms of the Brown representability theorem. 

Let $X$ be a non\-singular projective variety of dimension $n$ and write
$f\colon X\to Y=\Spec(\C)$ for the projection to a point.
In this case $f^!(\Oh_Y)=\om_X[n]$.
The above statement of Grothendieck duality implies that the functor
\begin{equation}
\label{sefu}
\Se_X(\blank)=(\blank\tensor\om_X)[n]
\end{equation}
is a Serre functor on $\D(X)$.

\subsection{Duality for quasi\-projective schemes}
\label{qpdual}

In order to apply Grothendieck duality on quasi\-projective schemes, we
need to restrict attention to sheaves with compact support. The {\em
support} of an object $E\in\D(X)$ is the locus of $X$ where $E$ is not
exact, that is, the union of the supports of the cohomology sheaves of
$E$. It is always a closed subset of $X$.

Given a scheme $X$, define the category $\Dc(X)$ to be the full
subcategory of $\D(X)$ consisting of complexes whose support is proper. 
Note that when $X$ itself is proper, $\Dc(X)$ is just the usual derived
category $\D(X)$.

If $f\colon X\to Y$ is a morphism of schemes of
finite Tor dimension, but not necessarily proper, 
then~(\ref{amos}) still holds for all objects in $\Dc(Y)$. Using this
we see that if $X$ is a non\-singular quasi\-projective variety of
dimension $n$, the category $\Dc(X)$ has a Serre functor given 
by~(\ref{sefu}).

\subsection{Crepant resolutions}
Let $X$ be a variety and $f\colon Y\to X$ a resolution of singularities.
Suppose that $X$ has rational singularities, that is, $f_*\Oh_Y=\Oh_X$ and
 \[
 \R^i f_* \Oh_Y=0 \qquad\text{for all } i>0.
 \]
Given a point $x\in X$ define $\D_x(Y)$ to be the full subcategory of $\Dc(Y)$
consisting of objects whose support is contained in the fibre $f\1(x)$.
We have the following categorical criterion for $f$ to be crepant.

\begin{lemma}
\label{crepancy}
If\/ $\D_x(Y)$ has trivial Serre functor for each $x\in X$,
then $X$ is Gorenstein and\/ $f\colon Y\to X$ is a crepant resolution.
\end{lemma}

\begin{pf}
The Serre functor on $\D_x(Y)$ is the restriction of the Serre
functor on $\Dc(Y)$.
Hence, by Section~\ref{qpdual}, the condition implies that for 
each $x\in X$ the restriction of the functor
$(\blank\tensor\om_Y)$ to the category $\D_x(Y)$ is isomorphic to the
identity. Since $\D_x(Y)$ contains the structure sheaves of all
fattened neighbourhoods of the fibre $f\1(x)$ this implies that the
restriction of $\om_Y$ to each formal fibre of $f$ is trivial.
To get the result, we must show that $\om_X$ is a line bundle
and that $f^*\om_X=\om_Y$.
Since $\om_X=f_*\om_Y$, this is achieved by the following lemma.
\end{pf}

\begin{lemma}
A line bundle $L$ on $Y$ is the pullback $f^*M$ 
of some line bundle $M$ on $X$ if
and only if the restriction of $L$ to each formal fibre of $f$ is trivial. 
Moreover, when this holds, $M=f_*L$.
\end{lemma}

\begin{pf}
For each point $x\in X$, the formal fibre of $f$ over $x$ is the fibre product
 \[
 Y\times_X \Spec(\widehat{\Oh}_{X,x}).
 \]
The restriction of the pullback of a line bundle from $X$ to each of these 
schemes is trivial because a line bundle has trivial formal stalks at points.

For the converse suppose that the restriction of $L$ to each of these formal
fibres is trivial. The theorem on formal functions shows that the completion
of the stalks of the sheaves $\R^i f_*\Oh_Y$ and $\R^i f_*L$ at any point
$x\in X$ are isomorphic for each $i$. Since $X$ has rational singularities
it follows that $\R^i f_*L=0$ for all $i>0$, and $M=f_*L$ is a line bundle
on $X$.

Since $f^*M$ is torsion free, the natural adjunction map 
$\eta\colon f^* f_*L \to L$ is injective, so there is a short exact sequence
\begin{equation}
\label{seseta}
0\to f^* f_*L \lRa{\eta} L\to Q\to 0.
\end{equation}
By the projection formula and the fact that $X$ is rational,
 \[
 \R^i f_*(f^* M)=M\tensor \R^i f_*\Oh_Y=0 \quad\mbox{for all }i>0.
 \]
The fact that $\eta$ is the unit of the adjunction for $(f^*,f_*)$ 
implies that $f_*\eta$ has a left inverse, and in particular is surjective.
Applying $f_*$ to~(\ref{seseta}) we conclude that $f_* Q=0$.

Using the theorem on formal functions again, we can deduce that 
 \[
 f_*(Q\tensor L\1)=0.
 \]
In particular, $Q\tensor L\1$ has no global sections. 
Tensoring~(\ref{seseta}) with $L\1$ gives a contradiction unless $Q=0$. Hence
$\eta$ is an isomorphism and we are done.
\end{pf}

%=======================================================================
\section{$G$-sheaves}
%=======================================================================

Throughout this section $G$ is a finite group acting 
on a scheme $X$ (on the left) by auto\-morphisms.
As in the last section, all schemes are of
finite type over $\C$. We list some results we need concerning the
category of sheaves on $X$ equipped with a compatible $G$ action,
or `$G$-sheaves' for short.
Since $G$ is finite, most of the
proofs are trivial and are left to the reader. The main point is that
natural constructions involving sheaves on $X$ are canonical, so commute
with automorphisms of $X$. 

\subsection{Sheaves and functors}

A $G$-sheaf $E$ on $X$ is a quasi\-coherent sheaf of
$\Oh_X$-modules together with a lift of the $G$ action to $E$.
More precisely, for each $g\in G$, there is a lift
$\lift{g}{E}\colon E\to g^*E$ satisfying $\lift{1}{E}=\id_E$
and $\lift{hg}{E}=g^*\left(\lift{h}{E}\right)\comp\lift{g}{E}$.

If $E$ and $F$ are $G$-sheaves, then there is a (right)
action of $G$ on $\Hom_X(E,F)$ given by
$\theta^g=\left(\lift{g}{F}\right)\1\comp g^*\theta\comp\lift{g}{E}$
and the space $\GHom_X(E,F)$ of $G$-invariant maps give the morphisms
in the Abelian categories $\Qco^G(X)$ and $\Coh^G(X)$ of
$G$-sheaves.

The category $\Qco^G(X)$ has enough injectives (Grothendieck
\cite[Proposition 5.1.2]{Gr}) so we may take $G$-equivariant injective
resolutions. Since $G$ is finite, if $X$ is a quasi\-projective scheme
there is an ample invertible $G$-sheaf on $X$ and so we may also
take \hbox{$G$-equivariant} locally free resolutions. The functors
$\GExt^i_X(\blank,\blank)$ are the $G$-invariant parts of
$\Ext^i_X(\blank,\blank)$ and are the derived functors of
$\GHom_X(\blank,\blank)$. Thus if $X$ is non\-singular of dimension $n$,
so that $\Qco(X)$ has global dimension $n$, then the category $\Qco^G(X)$
also has global dimension $n$.

The local functors $\sHom$ and $\tensor$ are defined in the obvious way on
$\Qco^G(X)$, as are pullback $f^*$ and pushforward $f_*$ for any
$G$-equivariant morphism of schemes $f\colon X\to Y$. Thus, for example,
$\lift{g}{f^*E}=f^*\lift{g}{E}$. Natural isomorphisms such as 
$\Hom_X(f^*E,F)\iso\Hom_Y(E,f_*F)$ are canonical, that is, commute with
isomorphisms of the base, and hence are $G$-equivariant. Therefore they
restrict to natural iso\-morphisms
 \[
 \GHom_X(f^*E,F)\iso\GHom_Y(E,f_*F).
 \]
In other words, $f^*$ and $f_*$ are also adjoint functors between the
categories $\Qco^G(X)$ and $\Qco^G(Y)$.

Similarly, the natural isomorphisms implicit in the projection formula,
flat base change, etc. are canonical and hence $G$-equivariant.

It seems worthwhile to single out the following point:

\begin{lemma}
\label{last}
Let $E$ and $F$ be $G$-sheaves on $X$.
Then, as a representation of $G$, we have a direct sum decomposition
 \[
 \Hom_X(E,F)=\bigoplus_{i=0}^k\GHom_X(E\tensor\rho_i,F)\tensor\rho_i
 \]
over the irreducible representations $\{\rho_0,\cdots,\rho_k\}$.
\end{lemma}

\begin{pf}
The result amounts to showing that 
\[
\GHom(\rho_i, \Hom_X(E,F))=\GHom_X(E\tensor\rho_i,F).
\]
Let $f\colon X\to Y=\Spec(\C)$ be projection to a point,
with $G$ acting trivially on $Y$ so that the map is equivariant.
Then $\Qco^G(Y)$ is just the category of $\C[G]$-modules.
Note that $\Hom_X(E,F)= f_*\sHom _{\Oh_X}(E,F)$ and
$f^*\rho_i=\Oh_X\tensor\rho_i$,
so that the adjunction between $f^*$ and $f_*$ gives
\begin{eqnarray*}
 \GHom_Y(\rho_i, f_* \sHom_{\Oh_X}(E,F))&=&
\GHom_X(\Oh_X\tensor\rho_i, \sHom_{\Oh_X}(E,F)) \\ 
&=&\GHom_X(E\tensor\rho_i,F),
\end{eqnarray*}
as required.
\end{pf}

\subsection{Trivial actions}
\label{sec:triv}

If the group $G$ acts trivially on $X$, then any $G$-sheaf $E$
decomposes as a direct sum
 \[
 E=\bigoplus_i E_i\tensor \rho_i
 \]
over the irreducible representations $\{\rho_0,\rho_1,\dots,\rho_k\}$ of $G$
(where $\rho_0=\mathbf1$ is the trivial representation). The sheaves $E_i$
are just ordinary sheaves on $X$. Furthermore,
$\GHom_X(E_i\tensor \rho_i,E_j\tensor \rho_j)=0$ for $i\ne j$. Thus the
category $\Qco^G(X)$ decomposes as a direct sum $\bigoplus_i\Qco^i(X)$ and
each summand is equivalent to $\Qco(X)$.

In particular, every $G$-sheaf $E$ has a fixed part $[E]^G$
and the functor 
 \[
 [\blank]^G\colon\Qco^G(X)\to\Qco(X)
 \]
is the left and right adjoint to the functor
 \[
 \blank\tensor\rho_0\colon\Qco(X)\to\Qco^G(X),
 \]
that is, `let $G$ act trivially'.
Both functors are exact.

\subsection{Derived categories}

The $G$-equivariant derived category $\D^G(X)$ is defined to be the 
full subcategory of the (unbounded) derived category of $\Qco^G(X)$
consisting of complexes with bounded and coherent cohomology.

The usual derived functors $\R\sHom$, $\Ltensor$, $\L f^*$ and $\R f_*$ may
be defined on the equivariant derived category, and, as for sheaves, the
standard properties of adjunctions, projection formula and flat base change
then hold because the implicit natural iso\-morphisms are sufficiently
canonical.

To obtain an equivariant Grothendieck duality we refer to Neeman's
results \cite{Ne}. Let $f\colon X\to Y$ be an equivariant morphism of
schemes. The only thing to check is that equivariant pushdown
$\R f_*$ commutes with small coproducts. This is proved exactly as in
\cite{Ne}. Then the functor $\R f_*$ has a right adjoint $f^!$, 
and~(\ref{amos}) holds when $f$ is proper and of finite Tor dimension.
Moreover the same result holds if $f$ is not proper, providing that we
restrict $\R f_*$ to the subcategories of objects with compact support.
Thus if $X$ is a non\-singular quasi\-projective variety of dimension $n$,
the full subcategory $\DcG(X)\subset\D^G(X)$ consisting of objects with
compact supports has the Serre functor
 \[
 \Se_X(\blank)\,=\,(\blank\tensor\om_X)[n],
 \]
where $\om_X$ is the canonical bundle of $X$ with its induced
$G$-structure.

\subsection{Indecomposability}
If $G$ acts trivially on $X$ then the results of Section~\ref{sec:triv} show
that $\D^G(X)$ decomposes as a direct sum of orthogonal subcategories indexed
by the irreducible representations of $G$. More generally it is easy to see
that $\D^G(X)$ is decomposable unless $G$ acts faithfully. We need the
converse of this statement.

\begin{lemma}
Suppose a finite group $G$ acts faithfully on a quasi\-projective variety $X$.
Then $\D^G(X)$ is indecomposable.
\end{lemma}

\begin{pf}
Suppose that $\D^G(X)$ decomposes as an orthogonal direct sum of two
subcategories $\A_1$ and $\A_2$. Any indecomposable object of $\D^G(X)$ lies in
either $\A_1$ or $\A_2$ and 
 \[
 \Hom_{\D^G(X)}(a_1,a_2)=0 \mbox{ for all } a_1\in \A_1, a_2\in\A_2. %%%
 \]
Since the action of $G$ is faithful, the general orbit is free. Let
$D=G\cdot x$ be a free orbit. Then $\Oh_D$ is indecomposable as a  
$G$-sheaf. Suppose without loss of generality that $\Oh_D$ lies in $\A_1$.

Let $\rho_i$ be an irreducible representation of $G$. The sheaf
$\Oh_X\tensor\rho_i$ is indecomposable in $\D^G(X)$ and there exists an
equivariant map $\Oh_X\tensor\rho_i\to\Oh_D$ so $\Oh_X\tensor\rho_i$ also
lies in $\A_1$. Any indecomposable $G$-sheaf $E$ supported in
dimension 0 has a section, so by Lemma~\ref{last} there is an equivariant
map $\Oh_X\tensor\rho_i\to E$, and thus $E$ lies in $\A_1$.

Finally given an indecomposable $G$-sheaf $F$, take an orbit
$G\cdot x$ contained in $\Supp(F)$ and let $i\colon G\cdot x\into X$
be the inclusion. Then $i_* i^*(F)$ is supported in dimension 0 and
there is an equivariant map $F\to i_* i^*(F)$, so $F$ also lies in
$\A_1$. Now $\A_2$ is orthogonal to all sheaves, hence is trivial.
\end{pf}

%=======================================================================
\section{The intersection theorem}
%=======================================================================

Our proof that $\GHilb{M}$ is non\-singular follows an idea developed in
Bridgeland and Maciocia \cite{Br2} for moduli spaces over K3 fibrations,
and uses the following famous and difficult result of commutative algebra:

\begin{thm}[Intersection theorem]
\label{com}
Let $(A,\m)$ be a local\/ \hbox{$\C$-algebra} of dimension $d$. Suppose that
 \[
 0\to M_s\to M_{s-1}\to\cdots\to M_0\to 0
 \]
is a nonexact complex of finitely generated free $A$-modules with each
homology module $H_i(M_{\blob})$ an $A$-module of finite length. Then $s\ge
d$. Moreover, if $s=d$ and $H_0(M_{\blob})\iso A/\m$, then
 \[
 H_i(M_{\blob})=0 \quad\text{for all $i\ne 0$,}
 \]
and $A$ is regular.
\end{thm}

The basic idea is as follows.
Serre's criterion states that any finite length $A$-module
has homological dimension $\geq d$ and that $A$ is regular
precisely if there
is a finite length $A$-module which has homological dimension
exactly $d$.
The intersection theorem gives corresponding statements
for complexes of $A$-modules with finite length homology.
As a rough slogan, ``regularity is a property of the derived category".
For the main part of the proof, see Roberts \cite{Ro1}, \cite{Ro2};
for the final clause, see \cite{Br2}.

We may rephrase the intersection theorem
 using the language of support and homological
dimension. If $X$ is a scheme and $E$ an object in $\D(X)$,
then it is easy to check
\cite{Br2} that, for any closed point $x\in X$,
 \[
 x\in\Supp E\iff \Hom_{\D(X)}^i(E,\Oh_x)\ne 0 \text{ for some $i\in\Z$.}
 \]
The {\em homological dimension} of a nonzero object $E\in\D(X)$, written
$\hd E$, is the smallest non\-negative integer $s$ such that $E$ is
isomorphic in $\D(X)$ to a complex of locally free sheaves on $X$ of
length $s$. If no such integer exists we put $\hd E =\infty$. One can
prove \cite{Br2} that if $X$ is quasiprojective, and $n$ is a nonnegative
integer, then $\hd E \le n$ if and only if there is an integer $j$ such
that for any point $x\in X$
 \[
 \Hom^i_{\D(X)}(E,\Oh_x)=0 \mbox{ unless }j\le i\le j+n.
 \]
The two parts of Theorem~\ref{com} now become the following (cf.\ 
\cite{Br2}).

\begin{cor}
\label{ab}
Let $X$ be a scheme and $E$ a nonzero object of\/ $\D(X)$. Then
 \[
 \codim(\Supp E)\,\le\,\hd E.
 \]
\end{cor}

\begin{cor}
\label{smooth}
Let $X$ be an irreducible $n$-dimensional scheme, and fix a point $x\in X$.
Suppose that there is an object $E$ of\/ $\D(X)$ such that for any point
$z\in X$, and any integer $i$,
 \[
 \Hom^i_{\D(X)}(E,\Oh_z)=0 \quad\mbox{unless }z=x \mbox{ and \/} 0\le i\le n.
 \]
Suppose also that $H_0(E)\iso\Oh_x$. Then $X$ is non\-singular at $x$
and $E\iso\Oh_x$.
\end{cor}

%=======================================================================
\section{The projective case}
%=======================================================================

\label{Sec6!Proj_case}

The aim of this section is to prove Theorem~\ref{second} under the
additional assumption that $M$ is projective. The quasi\-projective case
involves some further technical difficulties that we deal with in the
next section. Take notation as in the introduction.
We break the proof up
into 7 steps.

\Step1
Let $\pi_Y\colon Y\times M\to Y$ and $\pi_M\colon Y\times M\to M$ denote
the projections. The functor $\Phi$ may be rewritten
 \[
 \Phi(\blank)\iso \R\pi_{M_*}(\Oh_\sZ\tensor\pi_Y^*(\blank\tensor\rho_0)).
 \]
Note that $\Oh_\sZ$ has finite homological dimension, because $\sZ$ is flat
over $Y$ and $M$ is non\-singular. Hence the derived dual $\Oh_\sZ^{\dual}=
\R\sHom_{\Oh_{Y\times M}}(\Oh_\sZ,\Oh_{Y\times M})$ 
also has finite homological dimension
%(in fact $\Oh_\sZ^{\dual}= \om_{\sZ/Y}[-n]=p^!\Oh_Y$ by
%Grothendieck--Verdier duality),
and we may define another functor
$\Psi\colon\D^G(M)\to \D(Y)$, by the formula
 \[
 \Psi(\blank)=[\R \pi_{Y_*}(\sP\Ltensor \pi_M^*(\blank))]^G,
 \]
where $\sP=\Oh_\sZ^{\dual}\tensor\pi_M^*(\om_M)[n]$.

Now $\Psi$ is left adjoint to $\Phi$ because of the three standard adjunctions
described in Section~\ref{adjunctions}. The functor $\pi_M^*$ is the left
adjoint to $\R\pi_{M,^*}$. The functor $\blank\tensor\Oh_\sZ$ has the (left and
right) adjoint $\blank\tensor\Oh_\sZ^{\dual}$. Finally the functor $\pi_Y^!$
has the left adjoint $\R\pi_{Y_*}$ and
 \[
 \pi_Y^!(\blank)=\pi_Y^*(\blank)\tensor\pi_M^*(\om_M)[n].
 \]

\Step2
The composite functor $\Psi\comp\Phi$ is given by
 \[
 \R\pi_2{}_*(\sQ\Ltensor\pi_1^*(\blank)),
 \]
where $\pi_1$ and $\pi_2$ are the projections of $Y\times Y$ onto its
factors, and $\sQ$ is some object of $\D(Y\times Y)$. This is just
composition of correspondences  (see Mukai \cite[Proposition~1.3]{muk}).

If $i_y\colon\{y\}\times Y\into Y\times Y$ is the closed embedding then
$\L i_y^*(\sQ)=\Psi\Phi\Oh_y$, so that for any pair of points $y_1, y_2$,
 \begin{equation}
 \label{support}
 \Hom_{\D(Y\times Y)}^i(\sQ,\Oh_{(y_1,y_2)})=
 \Hom_{\D(Y)}^i(\Psi\Phi\Oh_{y_1},\Oh_{y_2}) =
 \GExt^i_M(\Oh_{Z_{y_1}},\Oh_{Z_{y_2}}),
 \end{equation}
using the adjunction for $(\Psi,\Phi)$. Our first objective is to show
that $\sQ$ is supported on the diagonal $\Delta\subset Y\times Y$,
or equivalently that the groups
in~(\ref{support}) vanish unless $y_1=y_2$. When $n=3$ this is the same
as the assumption (4.8) of Ito and Nakajima \cite{IN}.

\Step3
\label{bike}
Let $Z_1,Z_2\subset M$ be $G$-clusters. Then
 \[
 \GHom_M(\Oh_{Z_1},\Oh_{Z_2})
 =\begin{cases}
 \C &\text{if $Z_1=Z_2$,} \\
 0 &\text{otherwise.}
\end{cases}
 \]
To see this note that $\Oh_Z$ is generated as an $\Oh_M$ module
by any nonzero constant section.
But, since $H^0(\Oh_Z)$ is the regular representation of $G$,
the constant sections are precisely the \hbox{$G$-invariant} sections.
Hence any equivariant morphism maps a generator to a scalar multiple of
a generator and so is determined by that scalar.

Let $y_1$ and $y_2$ be distinct points of $Y$. Serre duality, together
with our assumption that $\om_M$ is locally trivial as a $G$-sheaf
implies that
 \[
 \GExt^n_M(\Oh_{Z_{y_1}},\Oh_{Z_{y_2}})=
 \GHom_M(\Oh_{Z_{y_2}},\Oh_{Z_{y_1}})=0,
 \]
so that
 \[
 \GExt^p_M(\Oh_{Z_{y_1}},\Oh_{Z_{y_2}})=0 \quad\text{unless $1\le p\le n-1$.}
 \]
Hence $\sQ$ restricted to $(Y\times Y)\setminus\Delta$ has homological
dimension $\le n-2$.

\Step4
Now we apply the intersection theorem.
If $y_1$ and $y_2$ are points of $Y$ such that $\tau(y_1)\ne \tau(y_2)$ then
the corresponding clusters $Z_{y_1}$ and $Z_{y_2}$ are disjoint, so that the
groups in~(\ref{support}) vanish. Thus the support of $\sQ|_{(Y\times
Y)\setminus\Delta}$ is contained in the subscheme $Y\times_X Y$. By
assumption this has codimension $>n-2$ so Corollary~\ref{ab} implies that
 \[
 \sQ|_{(Y\times Y)\setminus\Delta}\iso 0,
 \]
that is, $\sQ$ is supported on the diagonal.

\Step5
Fix a point $y\in Y$, and put $E=\Psi\Phi(\Oh_y)$. We proved above that
$E$ is supported at the point $y$. We claim that $H_0(E)=\Oh_y$. Note that
Corollary~\ref{smooth} then implies that $Y$ is non\-singular at $y$ and
$E\iso\Oh_y$.

To prove the claim, note that there is a unique map $E\to\Oh_y$, so we
obtain a triangle
 \[
 C\to E\to\Oh_y\to C[1]
 \]
for some object $C$ of $\D(Y)$. Using the adjoint pair
$(\Psi,\Phi)$, this gives a long exact sequence
 \begin{gather*}
 \cdots\to\Hom_{\D(Y)}^0(\Oh_y,\Oh_y)\to
 \Hom_{\D^G(M)}^0(\Phi\Oh_y,\Phi\Oh_y)\to \Hom_{\D(Y)}^0(C,\Oh_y) \\
 \to \Hom_{\D(Y)}^1(\Oh_y,\Oh_y)\lRa{\ep}
 \Hom_{\D^G(M)}^1(\Phi\Oh_y,\Phi\Oh_y)\to \cdots. 
 \end{gather*}
The homomorphism $\ep$ is just the Kodaira--Spencer map for the family of
clusters $\{\Oh_{Z_y}:y\in Y\}$ (Bridgeland \cite[Lemma~4.4]{Br1}), so is
injective. It follows that
 \[
 \Hom_{\D(Y)}^i(C,\Oh_y)=0\quad\text{for all $i\le0$.}
 \]
An easy spectral sequence argument (see \cite[Example~2.2]{Br1}), shows that
$H_i(C)=0$ for all $i\le0$. Taking homology sheaves of the above triangle
gives $H_0(E)=\Oh_y$, which proves the claim.

\Step6
We have now proved that $Y$ is non\-singular, and that for any pair of
points $y_1,y_2\in Y$, the homo\-morphisms
 \[
 \Phi\colon\Ext^i_Y(\Oh_{y_1},\Oh_{y_2})\to
 \GExt^i_M(\Oh_{Z_{y_1}},\Oh_{Z_{y_2}})
 \]
are isomorphisms. By assumption, the action of $G$ on $M$ is such that
$\om_M$ is trivial as a $G$-sheaf on an open neighbourhood
of each orbit $G\cdot x\subset M$. This implies that
 \[
 \Oh_{Z_y}\tensor\om_M\iso\Oh_{Z_y}
 \]
in $\Coh^G(M)$, for each $y\in Y$. Applying Theorem~\ref{tring2} shows
that $\Phi$ is an equivalence of categories.

\Step7
It remains to show that $\tau\colon Y\to X$ is crepant. Take a point
$x\in X=M/G$. The equivalence $\Phi$ restricts to give an equivalence
between the full subcategories $\D_x(Y)\subset\D(Y)$ and
$\D^G_x(M)\subset \D^G(M)$ consisting of objects supported on the
fibre $\tau\1(x)$ and the orbit $\pi\1(x)$ respectively.

The category $\D^G_x(M)$ has trivial Serre functor because $\om_M$ is
trivial as a $G$-sheaf on a neighbourhood of $\pi\1(x)$. Thus
$\D_x(Y)$ also has trivial Serre functor and Lemma~\ref{crepancy} gives
the result.

This completes the proof of Theorem~\ref{second} in the case that
$Y$ is projective.

%=======================================================================
\section{The quasi\-projective case}
%=======================================================================

In this section we complete the proof of Theorem~\ref{second}. Once again,
take notation as in the introduction. The only problem with the argument
of the last section is that when $M$ is not projective Grothendieck
duality in the form we need only applies to objects with compact support.
Thus if we restrict $\Phi$ to a functor
 \[
 \Phic\colon\Dc(Y)\to\DcG(M).
 \]
the argument of the last section carries through to show that $Y$ is
non\-singular and crepant and that $\Phic$ is an equivalence. It remains to
show that also $\Phi$ is an equivalence.

\Step8
The functor $\Phi$ has a right adjoint 
 \[
 \Ups(\blank)\,=\,[p_*\comp q^!(\blank)]^G
  \,=\,[\R\pi_{Y}{}_*(\om_{Z/M}\Ltensor\pi_M^*(\blank))]^G.
 \]
As before, the composition $\Ups\comp\Phi$ is given by
 \[
 \R\pi_2{}_*(\sQ\Ltensor\pi_1^*(\blank)),
 \]
where $\pi_1$ and $\pi_2$ are the projections of $Y\times Y$ onto its
factors, and $\sQ$ is some object of $\D(Y\times Y)$.

Since $\Phic$ is an equivalence, $\Ups\Phi\Oh_y=\Oh_y$ for any point
$y\in Y$, and it follows that $\sQ$ is actually the pushforward of a line
bundle $L$ on $Y$ to the diagonal in $Y\times Y$. The functor
$\Ups\comp\Phi$ is then just twisting by $L$, and to show that $\Phi$ is
fully faithful we must show that $L$ is trivial.

There is a morphism of functors $\ep\colon\id\to\Ups\comp\Phi$, which for
any point $y\in Y$ gives a commutative diagram
 \[
 \begin{CD}
\Oh_Y &@>\ep(\Oh_Y)>> &L \\
@VfVV && @VVL\tensor fV \\
\Oh_y &@>\ep(\Oh_y)>> &\,\Oh_y
 \end{CD}
 \]
where $f$ is non\-zero. Since $\ep$ is an isomorphism on the subcategory
$\Dc(Y)$, the maps $\ep(\Oh_y)$ are all isomorphisms, so the section
$\ep(\Oh_Y)$ is an isomorphism.

\Step9
The fact that $\Phi$ is an equivalence follows from Lemma~\ref{extra}
once we show that
 \[
 \Ups(E)\iso 0 \implies E\iso 0 \quad
 \text{for any object $E$ of $\D^G(M)$.}
 \]
Suppose $\Ups(E)\iso 0$. Using the adjunction for $(\Phi,\Ups)$,
 \[
 \Hom^i_{\D^G(M)}(B,E)=0 \quad\text{for all $i$,}
 \]
whenever $B\iso\Phi(A)$ for some object $A\in\D(Y)$. In particular, this holds
for any $B$ with compact support.

If $E$ is nonzero, let $D=G\cdot x$ be an orbit of $G$ contained in the
support of $E$. Let $i\colon D\into M$ denote the inclusion, a projective
equi\-variant morphism of schemes. Then the adjunction morphism
$i_* i^!(E)\to E$ is non\-zero, which gives a contradiction.

This completes the proof of Theorem~\ref{second}. \qed

%=======================================================================
\section{Nakamura's conjecture}
%=======================================================================

Recall that in Theorem~\ref{second} we took the space $Y$ to be an
irreducible component of $\GHilb{M}$. Note that when $Y$ is non\-singular
and $\Phi$ is an equivalence, $Y$ is actually a connected component. This
is simply because for any point $y\in Y$, the bijection
 \[
 \Phi\colon\Ext^1_Y(\Oh_y,\Oh_y)\to\GExt^1_M(\Oh_{Z_y},\Oh_{Z_y})
 \]
identifies the tangent space of $Y$ at $y$ with the tangent space of
$\GHilb{M}$ at $y$. In this section we wish to go further and prove that
when $M$ has dimension 3, $\GHilb{M}$ is in fact connected.

\paragraph{Proof of Nakamura's conjecture} Suppose by contradiction that
there exists a \hbox{$G$-cluster} $Z\subset M$ not contained among the
$\{Z_y: y\in Y\}$. Since $\Phi$ is an equivalence we can take
an object $E\in\Dc(Y)$ such that $\Phi(E)=\Oh_Z$. The argument of
Section~\ref{bike}, Step~3 shows that for any point $y\in Y$
 \[
 \Hom_{\D(Y)}^i(E,\Oh_y)=\GExt_M^i(\Oh_Z,\Oh_{Z_y})=0\quad
 \mbox{unless }1\le i\le 2.
 \]
This implies that $E$ has homological dimension 1, or more precisely, that
$E$ is quasi-isomorphic to a complex of locally free sheaves of the form
\begin{equation}
\label{ouch}
0\to L_2\lRa{f} L_1\to 0.
\end{equation}
But $\Oh_Z$ is supported on some $G$-orbit in $M$, so $E$ is supported on a
fibre of $Y$, and hence in codimension $\ge1$. It follows that the 
complex~(\ref{ouch}) is exact on the left, so $E\iso\coker f[1]$. 
In particular
$[E]=-[\coker f]$ in the Grothendieck group $K_{\mathrm c}(Y)$ of $\Dc(Y)$.

Let $y$ be a point of the fibre that is the support of $E$. By
Lemma~\ref{Kgrps} below, $[\Oh_{Z_y}]=[\Oh_Z]$ in $K^G_{\mathrm c}(M)$, so
that $[\Oh_y]=[E]$ in $K_{\mathrm c}(Y)$, since the equivalence $\Phi$
gives an isomorphism of Grothendieck groups.

Let $\Ybar$ be a non\-singular projective variety with an open inclusion
$i\colon Y\into \Ybar$. The functor $i_*\colon\Dc(Y)\to\D(\Ybar)$ induces a
map on K groups, so $[\coker f]=-[\Oh_y]$ in $K_{\mathrm c}(\Ybar)$. But
this contradicts Riemann--Roch, because if $L$ is a sufficiently ample line
bundle on $\Ybar$, then $\eu(\coker f\tensor L)$ and $\eu(\Oh_y\tensor L)$
are both positive.

\begin{lemma}
\label{Kgrps}
If $Z_1$ and $Z_2$ are two $G$-clusters on $M$ supported on the same orbit
then the corresponding elements $[\Oh_{Z_1}]$ and $[\Oh_{Z_2}]$ in the
Grothendieck group $K^G_{\mathrm c}(M)$ of $\D^G_{\mathrm c}(M)$ are equal.
\end{lemma}

\begin{pf}
We need to show that, as $G$-sheaves,
 $\Oh_{Z_1}$ and $\Oh_{Z_2}$ have composition series
with the same simple factors.
Suppose that they are both supported on the $G$-orbit $D=G\cdot x\subset M$
and let $H$ be the stabiliser subgroup of $x$ in $G$. 
The restriction functor
is an equivalence of categories 
from finite length $G$-sheaves supported on $D$
to finite length $H$-sheaves supported at $x$.
The reverse equivalence is the induction functor
$\left(\blank\otimes_{\C[H]}\C[G]\right)$.
Since the restriction of a $G$-cluster supported on $D$ is an
$H$-cluster supported at $x$,
it is sufficient to prove the result for 
$H$-clusters supported at $x$.

If $\{\rho_0,\cdots,\rho_k\}$ are the irreducible representations of $H$,
then we claim that the simple $H$-sheaves supported at $x$ are precisely
 \[
 \{S_i=\Oh_x\tensor\rho_i:0\le i\le k\}
 \]
These sheaves are certainly simple, 
since they are simple as $\C[H]$-modules.
On the other hand, any $H$-sheaf $E$ supported at $x$ has a 
nonzero ordinary sheaf morphism $\Oh_x\to E$.
By Lemma~\ref{last} there must be a nonzero $H$-sheaf morphism
$S_i\to E$, for some $i$, and, if $E$ were simple, then this would
have to be an isomorphism.

Thus a composition series as an $H$-sheaf is also a composition
series as a $\C[H]$-module. 
Hence all $H$-clusters supported at $x$ have the same composition
factors as $H$-sheaves, since as $\C[H]$-modules they are all
the regular representation of $H$.
\end{pf} 
 
%=======================================================================
\section{K theoretic consequences of equivalence}
%=======================================================================

In this section we put $M=\C^n$ and assume that the functor $\Phi$ is an
equivalence of categories. This is always the case when $n\le 3$. The main
point is that such an equivalence of derived categories immediately gives 
an isomorphism of the corresponding Grothendieck groups. 

\subsection{Restricting to the exceptional fibres}
Let  $\D^G_0(\C^n)$ denote the full subcategory of $\D^G(\C^n)$ consisting
of objects supported at the origin of $\C^n$. Similarly, let $\D_0(Y)$ denote the full subcategory of
$\D(Y)$ consisting of objects supported on the subscheme $\tau\1(\pi(0))$ of $Y$.

The equivalence $\Phi$ induces an equivalence
 \[
 \Phi_0\colon\D_0(Y)\to\D^G_0(\C^n),
 \]
so we obtain a diagram
 \[
 \renewcommand{\arraystretch}{1.5}
 \begin{array}{ccc}\D(Y) &\lRa{\Phi} &\D^G(\C^n) \\
 \bua && \bua \\
 \D_0(Y) &\lRa{\Phi} &\D^G_0(\C^n)
 \end{array}
 \]
in which the vertical arrows are embeddings of categories.

Note that the Euler characteristic gives natural bilinear forms on both
sides; if $E$ and $F$ are objects of $\D^G(\C^n)$ and $\D^G_0(\C^n)$
respectively, then we can compute
 \[
 \eu^G(E,F)=\sum_i(-1)^i \dim \Hom_{\D^G(\C^n)}(E,F[i]),
 \]
since the fact that the cohomology of $F$ has finite length
implies that the Ext groups are only nonvanishing in a finite interval.
Similarly, we can compute the ordinary Euler character on the left. The
fact that $\Phi$ is an equivalence of categories commuting with the shift
functors immediately gives
 \[
 \eu^G(\Phi(A),\Phi(B))=\eu(A,B),
 \]
for any objects $A$ of $\D(Y)$ and $B$ of $\D_0(Y)$.

\subsection{Equivalence of K groups}

Let $K(Y)$, $K^G(\C^n)$, $K_0(Y)$ and $K_0^G(\C^n)$ be the Grothendieck
groups of the corresponding derived categories. The equivalences of
categories from the last section immediately give isomorphisms of these
groups. The following lemma is proved in the same way as in
Gonzalez-Sprinberg and Verdier \cite[Proposition 1.4]{GV}.

\begin{lemma}
The maps that send a representation $\rho$ of $G$ to the $G$-sheaves
$\rho\tensor\Oh_{\C^n}$ and $\rho\tensor\Oh_0$ on $\C^n$ give ring
isomorphisms of the representation ring $R(G)$ with $K^G(\C^n)$ and
$K_0^G(\C^n)$ respectively.
\end{lemma}

We obtain a diagram of groups
 \[
 \renewcommand{\arraystretch}{1.5}
 \begin{array}{ccc}
 K(Y) &\lRa{\fie} & R(G) \\
 \scriptstyle{i}\bua\hphantom{\scriptstyle{i}} &&
 \hphantom{\scriptstyle{j}}\bua\scriptstyle{j} \\
 K_0(Y) &\lRa{\fie} & R(G).
 \end{array}
 \]
in which the horizontal maps are isomorphisms but the vertical maps
are not. In fact, if $Q$ is the representation induced by the
inclusion $G\subset\SL(n,\C)$, then the map $j$ is multiplication by
 \[
 r=\sum_{i=0}^n (-1)^i \Lambda^i Q\in R(G).
 \]
This formula is obtained by considering a Koszul resolution of $\Oh_0$
on $M$, as in \cite[Proposition 1.4]{GV}.
For example, in the case $n=2$ one has $r=2-Q$.

The above bilinear forms descend to give pairings on the Grothendieck
groups. These forms are non\-degenerate because if
$\{\rho_0,\cdots,\rho_k\}$ are the irreducible representations of $G$ then
the corresponding bases
 \[
 \{\rho_i\tensor\Oh_{\C^n}\}_{i=0}^k\subset K^G(\C^n)
\quad\text{and}\quad \{\rho_i\tensor\Oh_{0}\}_{i=0}^k\subset K^G_0(\C^n)
 \]
are dual with respect to the pairing $\eu^G(\blank,\blank)$. Applying
$\fie\1$ gives dual bases 
 \[
 \{\sR_i\}_{i=0}^k\subset K(Y)
\quad\text{and}\quad \{\sS_i\}_{i=0}^k\subset K_0(Y)
 \]
as in Ito and Nakajima \cite{IN}.

%=======================================================================
\section{Topological K~theory and physics}
%=======================================================================

With notation as in the introduction, suppose that $M$ is projective, and
further that $Y$ is non\-singular and $\Phi\colon\D(Y)\to\D^G(M)$ is an
equivalence. For example suppose that $n=2$ or $3$.

\subsection{K~theory and the orbifold Euler number}
Let $\K^*(Y)$ denote the topological complex K~theory of $Y$ and
$\K_G^*(M)$ the $G$-equi\-variant topological K~theory of $M$. There are
natural forgetful maps 
 \[
 \alpha_Y\colon K(Y)\to \K^0(Y)
 \quad\text{and}\quad \alpha_M\colon K^G(M)\to \K_G^0(M).
 \]
Since $\Phi$ and its inverse $\Psi$ are defined as correspondences,
we may define correspondences
 \[
 \fie\colon \K^*(Y)\to \K_G^*(M)
 \quad\text{and}\quad
  \psi\colon \K_G^*(M)\to\K^*(Y)
 \]
compatible with the maps $\alpha$, using the functors $\tensor$, $f^*$ and $f_*$
(also written $f_!$) on topological K~theory, which extend to equi\-variant 
K~theory, as usual, because they are canonical. Note that the definition and
compatibility of $f_*$ is non\-trivial; see \cite{AH} for more details. But
now the fact that $\Phi$ and $\Psi$ are mutually inverse implies that $\fie$
and $\psi$ are mutually inverse, that is, we have a graded isomorphism
\begin{equation}
\label{Kiso}
 \K^*(Y)\iso \K_G^*(M)
\end{equation}

Atiyah and Segal \cite{AS} observed that 
the physicists' orbifold Euler number of $M/G$ is 
the Euler characteristic of $\K_G^*(M)$, that is, 
 \[
 e(M,G)=\dim \K_G^0(M)\tensor\QQ -\dim \K_G^1(M)\tensor\QQ.
 \] 
On the other hand, since the Chern character gives a $\ZZ/2$ graded
isomorphism $\K^*(Y)\tensor\QQ \iso H^*(Y,\QQ)$, the Euler characteristic
of $\K^*(Y)$ is just the ordinary Euler number $e(Y)$ of $Y$.
Hence the isomorphism~(\ref{Kiso}) on topological K~theory provides a natural
explanation for the physicists' Euler number conjecture
 \[
 e(M,G)=e(Y).
 \]
This was verified in the case $n=2$ as a consequence of the original McKay
correspondence (cf.\ \cite{AS}). It was proved in the case $n=3$ by Roan
\cite{Roa} in the more general case of quasi\-projective Gorenstein
orbifolds, since the numerical statement reduces to the local linear
case $M=\C^3$, $G\subset\SL(3,\C)$.

\subsection{An example: the Kummer surface}

One of the first interesting cases of the isomorphism~(\ref{Kiso}) is when
$M$ is an Abelian surface (topologically, a 4-torus $T^4$), $G=\ZZ/2$
acting by the involution $-1$ and $Y$ is a K3 surface. In this case $Y$ is
a non\-singular Kummer surface, having 16 disjoint $-2$-curves
$C_1,\dots,C_{16}$ coming from resolving the images in $M/G$ of the 16
$G$-fixed points $x_1,\dots,x_{16}$ in $M$. Write $V=\{x_1,\dots,x_{16}\}$
for this fixed point set.

On the Abelian surface $M$ there are 32 flat line $G$-bundles,
arising from a choice of 2 $G$-actions on each of the 16 square
roots of $\Oh_M$. Each such flat line $G$-bundle $L(\rho)$ is
characterised by a map $\rho\colon V\to \FF_2=\{0,1\}$
such that at a fixed point $x\in V$ the group $G$ acts on 
the fibre $L_x$ with weight $(-1)^{\rho(x)}$. 
Now the set $V$ naturally has the structure of an affine
$4$-space over $\FF_2$ and the maps $\rho$ that occur are precisely the
affine linear maps, including the two constant maps corresponding to the
two actions on $\Oh_M$.

On the other hand, on the K3 surface $Y$ one may consider the lattice
$\ZZ^V\subset H^2(Y,\ZZ)$ spanned by $C_1,\dots,C_{16}$ and the smallest
primitive sublattice $\Lambda$ containing $\ZZ^V$. The elements of
$\Lambda$ give precisely the rational linear combinations of the divisors
$C_1,\dots,C_{16}$ which are themselves divisors. It is easy to see that
$\ZZ^V\subset\Lambda\subset(\frac12\ZZ)^V$ and it can also be shown that the
image of $\Lambda$ in the quotient $(\frac12\ZZ)^V/\ZZ^V\iso\FF_2^V$
consists of precisely the affine linear maps on $V$ (see Barth, Peters and
Van de Ven \cite[Chapter~VIII, Proposition 5.5]{BPV}).

We claim that under the correspondence $\Psi$,
the flat line $G$-bundle $L(\rho)$ is taken to the line
bundle $\Oh_Y(D(\rho))$, where
 \[
 D(\rho)=\frac12\Bigl( \sum_i \rho(x_i) C_i \Bigr).
 \]
To check the claim note that $\Oh_M$ is taken to $\Oh_Y$, and that, in the
local linear McKay correspondence for $\C^2/(\ZZ/2)$, the irreducible
representation of weight $-1$ is taken to the line bundle $\Oh(\frac12 C)$,
dual to the $-2$-curve $C$ resolving the singularity.
 
%=======================================================================

%=======================================================================

\bigskip \noindent
Tom Bridgeland, \\
Department of Mathematics and Statistics, \\
University of Edinburgh, King's Buildings, \\
Mayfield Road, Edinburgh EH9 3JZ, UK \\
e-mail: tab@maths.ed.ac.uk

\bigskip\noindent
Alastair King, \\
Department of Mathematical Sciences, \\
University of Bath, \\
Bath BA2 7AY, England \\
e-mail: a.d.king@maths.bath.ac.uk
%web: www.maths.bath.ac.uk/$\sim$masadk

\bigskip\noindent
Miles Reid, \\
Math Inst., Univ.\ of Warwick, \\
Coventry CV4 7AL, England \\
e-mail: miles@maths.warwick.ac.uk
%web: www.maths.warwick.ac.uk/$\sim$miles

%=======================================================================
\end{document}